\documentclass[10pt]{article}
\makeatletter

\@addtoreset{equation}{section} \makeatother
\newcommand{\R}{\mbox{I\hspace{-.25em}R}}

\def\R{\mathop{\rm I\kern  -,057cm R}}
\def\rn{\R^N}
\catcode`\@=11

\def\eqalign#1{\null\,\vcenter{\openup1\jot \m@th
  \ialign{\strut\hfil$\displaystyle{##}$&$\displaystyle{{}##}$\hfil
     &&\strut$\displaystyle{##}$\hfil&$\displaystyle{{}##}$
     \hfil\crcr#1\crcr}}\,}

\def\elle#1{L^{#1}(\Omega)}

\def\de{\delta}

\def\ct{{\cal T}}
\def\sob#1{W^{1,{#1}}_0(\Omega)}
\def\dive{{\rm div}}

\def\frac#1#2{{#1\over  #2}}
\def\vep{\varepsilon}
\def\r{{\bf R}}
\def\acca{H^1_0(\Omega)}
\def\dacca{H^{-1}(\Omega)}

\def\rn{{\bf R}^N}

\def\into{\int_\Omega}
\def\vfi{\varphi}
\def\qed{{\unskip\nobreak\hfil\penalty50
          \hskip2em\hbox{}\nobreak\hfil\mbox{\rule{1ex}{1ex} \qquad}
   \parfillskip=0pt
   \finalhyphendemerits=0\par\medskip}}
 \def\la{\lambda}
\def\be{\begin{equation}}
\def\ee{\end{equation}}
\def\rife#1{(\ref{#1})}
\def\proof{\noindent{\bf Proof.}\quad}

\newtheorem{theo}{Theorem}[section]
\newtheorem{rem}{Remark}[section]
\newtheorem{lem}{Lemma}[section]
\newtheorem{corol}{Corollary}[section]
\newtheorem{defi}{Definition}[section]

\begin{document}

\title{
  Uniqueness for unbounded solutions to stationary viscous Hamilton--Jacobi
equations\\ }
\author{ {\bf\large Guy Barles}, 
{\bf \large Alessio  Porretta}}
\vskip2em 

\date{}

\maketitle

\begin{small}
\noindent{\bf Abstract}  
We consider  a class of stationary viscous Hamilton--Jacobi equations as 
$$
\left\{\begin{array}{l}
\la\,u-{\rm div}(A(x) \nabla u)=H(x,\nabla u)\mbox{ in
}\Omega ,\\ u=0\mbox{ on }\partial\Omega\end{array}
 \right. 
$$
where $\la\geq 0$, $A(x)$ is a bounded and uniformly elliptic matrix and $H(x,\xi)$ is convex in $\xi$ and grows at most like 
$|\xi|^q+f(x)$, with $1<q<2$ and $f\in \elle {\frac N{q'}}$. Under such growth conditions solutions  are in general unbounded, and  there is not uniqueness of usual weak  solutions.  We prove  that   uniqueness holds in the restricted  class of  solutions satisfying  a suitable energy--type estimate, i.e. $(1+|u|)^{\bar q-1}\,u\in \acca$, for a  certain (optimal) exponent $\bar q$. This completes  the recent results in \cite{GMP}, where the existence of at least one  solution in this class has been proved. 
\end{small}
\vskip1em
\noindent {\bf  MSC}:{
 {\bf 35J60} (35R05, 35Dxx)}.
 \vskip1em
 \noindent {\it Running head: Viscous Hamilton--Jacobi equations} 
 \\
\section{Introduction}
In this paper we consider a class of 
elliptic  equations in a  bounded domain $\Omega \subset \R^{N}$, $N > 2$
\be\label{oo}
\left\{\begin{array}{l}
\lambda
u-\mbox{div}(A(x)\nabla u)=H(x,\nabla u)\mbox{ in
}\Omega ,\\ u=0\mbox{ on }\partial\Omega\end{array}
 \right. 
\ee
where the function $H(x,\xi)$ is convex and superlinear with respect to $\xi$. 

Equations of this type are sometimes referred to as stationary viscous
Hamilton--Jacobi equations and appear  in connection to  stochastic optimal control
problems. In that context, the convexity of $H$ is a natural assumption.

The model example which we are going to treat is the following
\begin{equation}\label{pbmod}
\left\{\begin{array}{l}
\lambda
u-\mbox{div}(A(x)\nabla u)=\gamma\,|\nabla u|^q+f(x)\mbox{ in
}\Omega ,\\ u=0\mbox{ on }\partial\Omega\end{array}
 \right.\end{equation} 
where  $q>1$, $\la\geq 0$ and  $A(x)=(a_{i,j}(x))$ is
a  matrix of $\elle\infty$
functions $a_{i,j}(x)$ satisfying uniform
ellipticity and boundedness conditions
\be\label{matr}
\alpha |\xi|^2\leq A(x)\xi\cdot \xi\leq
\beta |\xi|^2\qquad \forall \xi\in
\rn\,,\quad {\rm a.e.}\quad x\in \Omega\,.
\ee
Without loss of generality, we let $\gamma >0$.  
We draw our attention to the 
\lq\lq subcritical\rq\rq\ case, namely  $q<2$, and, more precisely,  to the question of uniqueness  of {\sl  unbounded solutions}.

Let us first recall that  some regularity condition is needed
on $f$ in order that problem \rife{pbmod} admits a
solution. In the class of Lebesgue spaces, this
condition amounts to ask that 
\be\label{regf}
f\in \elle{\frac
N{q'}},
\ee
where $q'$ is the conjugate exponent of $q$, i.e. $\displaystyle q' = \frac{q}{q-1}$.
When $q<2$, \rife{regf} implies  that $f\in \elle m$ with $m<\frac N2$, hence solutions
are expected to be unbounded. Moreover since, by Sobolev embedding theorem, one has
\be\label{ac}
\elle{\frac N{q'}}\subset \dacca \quad \iff \quad q\geq 1+\frac 2N \, ,
\ee
the value $q=1+\frac 2N$ is  a critical one. Indeed, the solutions  
belong to $\acca$ only if $q\geq 1+\frac 2N$, when $q$ is below this value solutions are
not only unbounded but have not even finite energy and should be defined in a suitable
generalized sense. 

The fact that \rife{regf} is a necessary condition for having solutions can be easily justified by a heuristic argument: if $A(x)=I$, i.e. in case of the Laplace operator, the Calderon--Zygmund regularity  
implies that 
$$
-\Delta u\in \elle m\rightarrow u\in
W^{2,m}(\Omega)\rightarrow |\nabla u|\in
\elle{m^*}\, ,
$$
where $m^*$ is associated to $m$ through the Sobolev embedding, i.e., for $N>m$, $m^* = \displaystyle \frac{Nm}{N-m}$.
In order to be consistent with \rife{pbmod}, this
means that $f\in \elle m $ and $|\nabla u|\in \elle{qm}$ so that  one needs
$qm=m^*$, i.e.  $m=\frac N{q'}$. We refer the reader to \cite{AP}, \cite{HMV} for rigorous
and sharper necessary conditions  on
$f$  in order to have weak solutions.  It is important to   recall that if
$\la=0$ the data  $f$, $\gamma$, $\alpha$ must also satisfy a  size condition in order
that a solution exists.

Pioneering results for such kind of equations were given by P.L. Lions (\cite{Li1},
\cite{Li2}), mainly in case of Lipschitz solutions and including $q>2$. Existence results
for the case
$q=2$ can be found in several works, among which we recall the series of papers by L.
Boccardo, F. Murat, J.P. Puel (see e.g.  
\cite{BMP}, \cite{BMP1}) and more recently,   assuming  $f$ in $\elle{\frac N2}  $,  
in  \cite{FM}, \cite{DGP}.

Under assumption \rife{regf} with $1<q<2$, the
existence of a solution for problems as  \rife{pbmod} has been
recently proved in \cite{GMP} if either $\la>0$ or  $\la=0$ and  
a size condition  is satisfied
\be\label{size}
\gamma^{\frac1{q-1}}\,\|f\|_{\elle {\frac N{q'}}}<
\alpha^{q'}C_*\,,
\ee
where $C_*$ only depends on    $q$ and  $N$.

In this paper we deal with the problem  of 
uniqueness of solutions. 
Up to now, uniqueness  results for problems like \rife{pbmod}
have been proved in the Lipschitz  case (\cite{Li1}) and in  \cite{Bar-Mu},
\cite{Bar2}
 if either solutions are  bounded or
$q=2$. Note that these two cases share a common
feature, which is that   $f$ is required to be in
$\elle m$,
$m\geq
\frac N2$: for less summable $f$ as we consider,   the approach of these previous
papers   seems not to 
work.  Further results when $q
\leq 1+\frac 2N$ can be found in \cite{BMMP1}, \cite{BMMP2}.

When dealing with the question of uniqueness, 
one has to consider the following well--known counterexample  (see also  \cite{Li1},
\cite{AP}) for $q>\frac N{N-1}$
\be\label{count}
\eqalign{ &\hbox{ 
$u(x)=C_\alpha\,(|x|^\alpha-1)$ ($\alpha=-\frac{2-q}{q-1}$, 
$C_\alpha=\frac{(N+\alpha-2)^{\frac1{q-1}}}{|\alpha|}$) solves}
\cr
&
\cases{-\Delta u=|\nabla u|^q &in ${\cal D}'(B_1(0))$, \cr u\in \sob q&\cr}\cr}
\ee
This shows that   uniqueness
does not hold in the   class of weak solutions
$u$ in $\sob q$, and, if $q\geq 1+\frac 2N$, not even  in $\acca$ (one can check that $u\in \acca$
in this case). It is then natural to look for   a suitable   class of solutions in which 
problem \rife{pbmod} is well--posed.   
A  linearization argument would suggest that there is uniqueness in the class
\be\label{cla}
\hbox{$u$ sol. of \rife{pbmod}: $|\nabla u|\in
\elle{N(q-1)}$.}
\ee
On the other hand, if $q > 1+\frac 2N$ (which gives $N(q-1)>2$), the existence of 
such kind of  solutions
can not be obtained unless the Calderon--Zygmund regularity theorem applies; 
thus, in order to deal with    general (bounded measurable) coefficients
$a_{i,j}$, this approach is not reasonable. 

Our main purpose here is to prove the uniqueness of solutions of \rife{pbmod} in a
 regularity class 
which is consistent with the existence results  available from 
\cite{GMP}. In this latter paper  it has been proved  that   a  
natural class of solutions for which both {\sl a priori
estimates} and {\sl existence} hold is given through the extra energy  condition
\be\label{clac}
\hbox{$u$ sol. of \rife{pbmod}: $(1+|u|)^{\bar q-1}\,u\in \acca$, with $\bar q=
\frac{(N-2)(q-1)}{2(2-q)}$.}
\ee
We are going to prove that this regularity is precisely what is needed to  select a unique
solution, so that problem \rife{pbmod} is actually well--posed in this class. 

Our main result concerns the case  $q\geq 1+\frac 2N$,
which corresponds to $\acca$ solutions (see \rife{ac}).

\begin{theo}\label{1}
Let  $ 
1+\frac 2N\leq q<2$. Assume  \rife{matr}, 
\rife{regf} and that

(i) either $\la >0$

(ii) or $\la=0$  and \rife{size} is satisfied.
 
\noindent Then problem \rife{pbmod} has  only one (distributional) solution
$u$ such that
$(1+|u|)^{\bar q-1}u\in \acca$, with
$\bar q= \frac{(N-2)(q-1)}{2(2-q)}$. 
\end{theo}

Note that
the function $u$ in the counterexample
\rife{count} satisfies $(1+|u|)^{r-1}\,u\in \acca $ for any $r<\bar q$ but not for $r = \bar q$, which proves the
optimality of our result. Observe also 
that $\bar q$ tends to infinity as $q\to 2$, which is consistent with the case $f\in
\elle{\frac N2}$, for which existence and uniqueness have been  proved (see
\cite{FM},  \cite{Bar2} respectively) in the class
of solutions
$u$ such that $\exp(\mu\,u)-1\in \acca$ for  a
suitable constant $\mu$. 

We leave to Section 2 the proof of Theorem \ref{1}, actually in 
a generalized version which includes problem \rife{oo} where $H$ is convex and
satisfies similar growth conditions. Some extensions to Neumann boundary conditions as
well as to the case of unbounded domains is also discussed.

In Section 3 we deal with the case  
$\frac N{N-1}<q<  1+\frac 2N$, which corresponds to $f\in \elle m$
with $1<m<\frac{2N}{N+2}$. A similar result as Theorem \ref{1} is proved, but since, in
this case, solutions do not belong to
$\acca$, we use a slightly stronger  formulation than the distributional one, namely
uniqueness is proved for so--called renormalized solutions (still in the class
\rife{clac}). This notion (see Definition \ref{RE}), first introduced  in \cite{DL} for
transport equations, is now currently used in several different contexts when dealing with
solutions of infinite energy.  

Still in Section 3, we prove in fact a more general uniqueness result when $q$ is
below the critical value
$1+\frac 2N$. Indeed, we will see
that if
$q\leq 1+\frac 2N$  then 
the regularity
\rife{clac} implies \rife{cla}. This fact allows to prove  uniqueness   through
a  simpler linearization principle, which  does not need any convexity argument and
which can be applied to more general situations like, for instance, nonlinear operators
(see Theorem \ref{nonlin}).
Note that the limiting value $q=1+\frac 2N$ is also admitted here; actually, \rife{clac}
and \rife{cla} coincide in that case with $u\in \acca$. On the other hand, as
mentioned before, this argument was not possible for $q>1+\frac2 N$ since \rife{cla} will
no more be true in general.

Finally, some further   remarks will be discussed  at the end of Section 3,
including   the case 
$q<\frac N{N-1}$, where uniqueness holds simply in   $\sob q$.

\section{The case ${\bf q\geq 1+\frac2N}$: finite energy solutions}

\vskip1em

We consider  a natural generalization of
\rife{pbmod}, namely the following equation
\be\label{pb}
\left\{\begin{array}{l}
\lambda
u-\mbox{div}(A(x)\nabla u)=H(x,\nabla u) \mbox{ in }\Omega
,\\ u=0\mbox{ on }\partial\Omega\end{array}
 \right.
\ee
We still   assume that  $\la\geq 0$, that $A(x)$
satisfies \rife{matr} and that
$H(x,\xi)$ is a Carath\'eodory function
satisfying
 \be\label{h1}
\xi\mapsto H(x,\xi)\quad\hbox{is convex, for a.e.
$x\in \Omega$,}
 \ee
and   the growth condition
\be\label{h2}
\hbox{$\exists\,\, q\in (\frac N{N-1},2)$:}\quad |H(x,\xi)|\leq
\gamma\,|\xi|^q+f(x)\,,\quad
\gamma>0\,,\,\quad f(x)\in \elle{\frac N{q'}}.
\ee
Note that this assumptions include the possibility that the equation contains transport
terms; indeed, the basic choice for $H$ is  
$$
H(x,\nabla u)= b(x)\cdot \nabla
u + \gamma(x) |\nabla u|^p +f(x)
$$
where $b \in \left[\elle N\right]^N$, $\gamma \in \elle r$
and 
$f\in
\elle{\frac{Nr(p-1)}{pr-N}}$, with $r\in (N,+\infty]$ and $
\frac{N(r-1)}{r(N-1)}<p<2-\frac Nr$.
\vskip1em

In virtue of  \rife{h2} and \rife{ac}, assuming $q\geq 1+\frac 2N$ corresponds to having
data in
$\dacca$, so  that we can reasonably talk of $\acca$ weak  solutions.

\begin{defi}   We say that $u\in \acca$ is  a weak  subsolution  
of
\rife{pb} if   $H(x,\nabla u)\in \elle1$ and 
\be\label{WE}
\begin{array}{rl}
&
\la\,\into u\,\xi\,dx + \into A(x)\nabla u\nabla \xi\,dx \leq  \into H(x,\nabla
u)\xi\,dx \\ 
\noalign{\medskip}&\qquad
\qquad \forall \xi \in \acca\cap \elle\infty\,,\quad \xi\geq 0\,.
\end{array}
\ee
A super-solution 
of
\rife{pb} is defined if the opposite inequality holds. A function $u$ being both a  sub and
a  super-solution is said to be a weak solution of
\rife{pb}. 
\end{defi}

Our proof of the comparison principle for sub and super-solutions of
\rife{pb}   relies on two basic ideas: the first one
is that if 
\be\label{ineq}
\cases{
-\dive(A(x)\nabla w)\leq |\nabla w|^q&\cr (w^+)^{\bar
q}\in \acca\,,\cr}
\ee
then $w\leq 0$; in other words, the homogeneous problem  
has only the trivial solution in this regularity
class. Secondly, we aim at applying inequality
\rife{ineq} to
(a small perturbation of)
 the difference of two solutions $u -v$. In
order to obtain this inequation, we use a   
convexity argument, which gives account of assumption
\rife{h1}.  A further technical tool will be required
in order to justify some  regularity  claimed on
$u-v$: here we apply  a truncation argument.

In order to do that we introduce the following truncation function
\be\label{deft}
T_n(s)=\int_0^s \theta_n(\xi)d\xi\,,\quad
\theta_n(\xi)=\cases{1&if $|\xi|<n$\cr
\frac{2n-|\xi|}n  &if $n<|\xi|<2n$,\cr 0&if
$|\xi|>2n$\cr}
\ee
and we start by giving a sort of {\it  renormalization
principle}  for the  \lq\lq truncated\rq\rq\ equation.

\begin{lem}\label{ren}  Let $u\in \acca$ be a weak subsolution 
of \rife{pb}. Then
$u$ satisfies, for any nonnegative $\xi\in \acca\cap
\elle\infty$ and for every $n$
\be\label{req}
\eqalign{ &
\la\into
T_n(u)\xi\,dx+\into A(x)\nabla
T_n(u)\nabla \xi \,dx\cr  &\leq \into H(x,\nabla 
T_n(u)) \xi\,dx+
\langle I_{n},\xi\rangle\,,\cr}
\ee
where $I_{n}$ is  defined as
\be\label{mun}
\eqalign{ &
\langle
I_{ n},\xi\rangle=
\frac1n \int\limits_{\{n<u<2n\}}\!\!\!\!
A(x)\nabla u\nabla u\,\xi\,dx -
\frac1n\int\limits_{\{-2n<u<-n\}}\!\!\!\!A(x)\nabla
u\nabla u\,\xi\,dx \cr &
+\la\into (T_n(u)-u\theta_n(u))\xi\,dx+\into
(H(x,\nabla u)\theta_n(u)-H(x,\nabla
T_n(u)))\xi\,dx\,.\cr}
\ee
If moreover 
$|u|^{\bar q-1}u\in \acca$, where $\bar
q=\frac{(N-2)(q-1)}{2(2-q)}$, we have
\be\label{tail}
\eqalign{&
\lim\limits_{n\to +\infty} n^{2\bar q-1}\,\|
I_{n}\|_{L^1(\Omega)}=0\,.\cr}
\ee
\end{lem}

\proof Let $\xi\in \acca\cap \elle\infty$,
$\xi\geq 0$, and let $n>0$. Multiplying equation
\rife{pb}  by
$\theta_n(u)\xi$ we obtain 
\be\label{renp}
\eqalign{ &
\la\into u\,\theta_n(u)\xi\,dx+\into
A(x)\nabla u\nabla \xi \theta_n(u)\,dx =\into
H(x,\nabla  u)\theta_n(u) \xi\,dx\cr &\quad +
\frac1n \int\limits_{\{n<u<2n\}}\!\!\!\!
A(x)\nabla u\nabla u\,\xi\,dx -
\frac1n\int\limits_{\{-2n<u<-n\}}\!\!\!\!A(x)\nabla
u\nabla u\,\xi\,dx\cr}
\ee
Recalling that $\theta_n(u)=T_n'(u)$  and defining
$I_n$ as in \rife{mun} we have obtained \rife{req}.
Now let $u$ be such that $|u|^{\bar q-1}u\in \acca$, where $\bar
q=\frac{(N-2)(q-1)}{2(2-q)}$. We have, by definition
of $I_n$
\be\label{in}
\eqalign{ &\qquad\qquad
n^{2\bar q-1}\|I_n\|_{\elle1}\leq\cr & n^{2\bar
q-1}\left [\frac1n
\int\limits_{\{n<u<2n\}}\!\!\!\! A(x)\nabla u\nabla
u\,dx+
\frac1n\int\limits_{\{-2n<u<-n\}}\!\!\!\!A(x)\nabla
u\nabla u\,dx\right]
\cr &
+\la\,n^{2\bar q-1}\!\into \!\!|T_n(u)-u\theta_n(u) |dx+
n^{2\bar q-1}\!\into \!\! |H(x,\nabla
u)\theta_n(u)-H(x,\nabla T_n(u))|dx\cr}
\ee
Observe that 
\be\label{in1}
\eqalign{ &
n^{2\bar
q-1}\,\frac1n
\int\limits_{\{n<u<2n\}}\!\!\!\!\!\!\!\!  A(x)\nabla u\nabla
u\,dx\leq \int\limits_{\{n<u<2n\}}\!\!\!\!\!\!\!\!
A(x)\nabla u\nabla u\,|u|^{2\bar q-2}dx\cr &\quad \qquad \qquad 
\leq \frac{\beta}{{\bar q}^2}\int\limits_{\{n<u<2n\}}
\!\!\!\! |\nabla |u|^{\bar
q}|^2\,dx\mathop{\to}^{n\to +\infty}0\,,\cr}
\ee
and similarly 
\be\label{in2}
n^{2\bar
q-2}\int\limits_{\{-2n<u<-n\}}\!\!\!\!A(x)\nabla
u\nabla u\,dx\mathop{\to}^{n\to +\infty}0
\ee
We also have, using  Young's inequality, and by
definition of $\bar q$, 
$$
\eqalign{ &
|\nabla u|^q\,|u|^{2\bar q-1}\leq \frac q2\,|\nabla
u|^2|u|^{2(\bar q-1)}+\frac{2-q}2 |u|^{2\bar
q+\frac{2(q-1)}{2-q}}\cr &
=  \frac q{2{\bar q}^2}\,|\nabla
(|u|^{\bar q-1}u)|^2 +\frac{2-q}2 (|u|^{\bar
q})^{2^*}\,.\cr}
$$
Since $|u|^{\bar q-1}u\in
\acca$, we conclude that 
\be\label{gra}
|\nabla u|^q\,|u|^{2\bar q-1}\in \elle1
\ee
Similarly, since  $(2\bar q-1)(\frac
N{q'})'=2^*\,\bar q$ we have that $|u|^{2\bar
q-1}\in \elle{(\frac N{q'})'}$; since $f\in
\elle{\frac N{q'}}$ we deduce that 
\be\label{etf}
f\, |u|^{2\bar
q-1}\in\elle1
\ee
Now we have, by definition of $\theta_n$ and
$T_n$,
$$
|H(x,\nabla
u)\theta_n(u)-H(x,\nabla T_n(u))|\leq
[|H(x,\nabla u)| +|H(x,\nabla
T_n(u))|]\chi_{\{|u|>n\}}\,,
$$
hence, using the growth assumption \rife{h2},
$$
n^{2\bar q-1}|H(x,\nabla
u)\theta_n(u)-H(x,\nabla T_n(u))|\leq
2 (\gamma|\nabla u|^q+f(x)) |u|^{2\bar
q-1}\chi_{\{|u|>n\}}\,.
$$
Thanks to \rife{gra}--\rife{etf} we conclude
\be\label{in3}
n^{2\bar q-1}\into  |H(x,\nabla
u)\theta_n(u)-H(x,\nabla
T_n(u))|\,dx\mathop{\to}^{n\to +\infty}0\,.
\ee
Finally, since $u\in \elle{2^*{\bar q}}$,
\be\label{in4}
n^{2\bar q-1}\into |T_n(u)-u\theta_n(u) |\,dx\leq 
2\into |u|^{2\bar q}\chi_{\{|u|>n\}}\,dx
\mathop{\to}^{n\to +\infty}0.
\ee
From \rife{in}, \rife{in1}, \rife{in2}, \rife{in3},
\rife{in4} we
 get \rife{tail}.
\qed

\begin{rem}\rm
Clearly if $u\in \acca$ is a
super-solution  of \rife{pb} then \rife{req} holds with the opposite sign.
In particular, if $u$ is a  weak solution of \rife{pb}, then
$$
\eqalign{ &\la\into T_n(u)\,\xi\,dx+
\into A(x)\nabla T_n(u)\nabla \xi \,dx\cr  &=
\into H(x,\nabla T_n(
u)) \xi\,dx+
\langle I_{n},\xi\rangle\,,\cr}
$$
with $n^{2\bar q-1}\|I_n\|_{\elle1}\to 0$ provided
$|u|^{\bar q-1}u\in \acca$.
\end{rem}

We come to the main comparison result.

\begin{theo}\label{comp}
Assume \rife{matr}, \rife{h1}, \rife{h2} with $q\geq 1+\frac 2N$. Let  $\la>0$.
If $u$
and $v$ are respectively a  subsolution and a  super-solution of \rife{pb} such that
$(1+|u|)^{\bar q-1}u\in \acca$ and $(1+|v|)^{\bar q-1}v\in \acca$, then we have $u\leq v$ in $\Omega$. 

In particular,  problem
\rife{pb} has a unique weak solution $u$ such that
$(1+|u|)^{\bar q-1}u \in
\acca$.
\end{theo}

\proof
From
Lemma \ref{ren} we obtain that
\be\label{ueq}
\eqalign{ &
\la\into T_n(u)\,\xi\,dx+\into
 A(x)\nabla T_n(u)\nabla \xi \,dx\leq\into
H(x,\nabla T_n(u))\xi\,dx + 
\langle I_{n}^{u},\xi\rangle\,,\cr}
\ee
and
\be\label{veq}
\eqalign{ &
\la\into T_n(v)\,\xi\,dx+\into
 A(x)\nabla T_n(v)\nabla \xi \,dx\geq\into
H(x,\nabla T_n(v))\xi\,dx+
\langle I_{n}^{v},\xi\rangle\,,\cr}
\ee
for every $\xi\in \acca\cap \elle\infty$, $\xi\geq 0$.

Let now $\vep>0$ be fixed. 
Subtracting  
\rife{ueq} and \rife{veq}, we obtain 
$$
\eqalign{ &
\la\into [T_n(u)-(1-\vep)T_n(v)]\,\xi\,dx+\into
 A(x)\nabla (T_n(u)-(1-\vep)T_n(v))\nabla \xi
\,dx\cr &\quad \leq\into [H(x,\nabla
T_n(u))-(1-\vep)H(x,\nabla
T_n(v))]\xi\,dx+
\langle I_{n}^{u},\xi\rangle-(1-\vep)\langle I_{n}^{v},\xi\rangle\,.\cr}
$$
Now we use the convexity assumption on $H$, which
gives
$$
H(x,p)\leq (1-\vep)H(x,\eta)+\vep
H(x,\frac{p-(1-\vep)\eta}\vep)\,,\quad \forall
p\,,\eta\in \rn\,.
$$
With $p=\nabla T_n(u)$, $\eta = \nabla T_n(v)$ we
obtain
$$
H(x,\nabla T_n(u))-(1-\vep)H(x,\nabla
T_n(v))\leq \vep H\left(x,\frac{\nabla
T_n(u)-(1-\vep)\nabla T_n(v)}\vep\right)\,,
$$
hence we have
\be\label{prew}
\eqalign{ &
\la\into [T_n(u)-(1-\vep)T_n(v)]\,\xi\,dx+\into
 A(x)\nabla [T_n(u)-(1-\vep)T_n(v)]\nabla \xi
\,dx\cr &\quad \leq \vep\into H\left(x,\frac{\nabla
[T_n(u)-(1-\vep)T_n(v)]}\vep\right)\xi\,dx+
\langle I_{n}^{u},\xi\rangle-(1-\vep)\langle I_{n}^{v},\xi\rangle\,.\cr}
\ee
We define now
\be\label{defw}
w_n=T_n(u)-(1-\vep)T_n(v)-\vep\vfi\,,
\ee
where $\vfi$ is  a positive function that belongs to 
$\acca$ and  will be chosen later.
From \rife{prew}  we obtain
\be\label{prefi}
\eqalign{ &
\la\into w_n\,\xi\,dx+\into
 A(x)\nabla w_n\nabla \xi
\,dx \leq \vep \into H\left(x,\frac{\nabla
w_n}\vep+\nabla \vfi\right)\xi\,dx\cr &\quad\quad
-\vep
\left[\into \la\vfi\,\xi+A(x)\nabla
\vfi\nabla \xi\,dx\right]+
\langle I_{n}^{u},\xi\rangle-(1-\vep)\langle I_{n}^{v},\xi\rangle\,.\cr}
\ee
Using assumption \rife{h2}   we
have
\be\label{de}
\eqalign{ &
H\left(x,\frac{\nabla
w_n}\vep+\nabla \vfi\right)\leq
\gamma\left|\frac{\nabla w_n}\vep+\nabla
\vfi\right|^q+f(x)\cr &\qquad \leq (\gamma+\de)
\left|\nabla
\vfi\right|^q+C_\de\left|\frac{\nabla
w_n}\vep\right|^q+f(x)\,,\cr}
\ee
where $\de$ is any positive constant.

Using \rife{de} in \rife{prefi} we obtain
\be\label{ec}
\eqalign{ &
\la\into w_n\,\xi \,dx+\into
 A(x)\nabla w_n\nabla \xi 
\,dx \leq \frac{C_\de}{\vep^{q-1}} \into
\left|\nabla w_n\right|^q\,\xi \,dx\cr
 &\quad\quad
-\vep
\left[\into \la\vfi\,\xi +A(x)\nabla
\vfi\nabla
\xi\,dx-(\gamma+\de)|\nabla
\vfi|^q\xi-f(x)\xi\right]\cr 
& \qquad \qquad +\langle I_{n}^{u},\xi\rangle-(1-\vep)\langle I_{n}^{v},\xi\rangle\,.\cr}
\ee
We choose now $\vfi$ as a solution of 
\be\label{aux}
\cases{\la\,\vfi-\dive(A(x)\nabla
\vfi)=(\gamma+\de)|\nabla \vfi|^q+f(x)& in
$\Omega$,\cr
(1+|\vfi|)^{\bar
q-1}\vfi\in
\acca.&\cr}
\ee
The existence of such a function $\vfi$ is proved in
\cite{GMP}. Moreover, we have that $\vfi\geq 0$ (since $f\geq 0$ from \rife{h2}).
Thanks to \rife{aux} we obtain from \rife{ec}
\be\label{ec1}
\eqalign{  
\la\into w_n\,\xi \,dx+\into
 A(x)\nabla w_n\nabla \xi 
\,dx \leq &\frac{C_\de}{\vep^{q-1}} \into
\left|\nabla w_n\right|^q\,\xi \,dx\cr &\qquad  +\langle
I_{n}^{u},\xi\rangle-(1-\vep)\langle I_{n}^{v},\xi\rangle\,.\cr}
\ee
For $l>0$, we choose  in \rife{ec1} $\xi=\xi_n$ defined as
$$
\xi_n= [(w_n-l)^+]^{2\bar q-1}
$$
Note that $\xi_n$ is a  positive function, and it
belongs to $\acca\cap \elle\infty$. Moreover, by the definition of $w_n$ in \rife{defw},
we have 
$$
\|\xi_n\|_{\elle\infty}\leq (2n)^{2\bar
q-1},
$$
so that we can apply Lemma \ref{ren} for $u$ and $v$
and get
$$
\eqalign{ &
|\langle I_{n}^{u},\xi_n\rangle|\leq (2n)^{2\bar
q-1}\|I_{n}^{u}\|_{L^1(\Omega)}\mathop{\to}^{n\to
+\infty}0\cr  & 
|\langle I_{n}^{v},\xi_n\rangle|\leq (2n)^{2\bar
q-1}\|I_{n}^{v}\|_{L^1(\Omega)}\mathop{\to}^{n\to
+\infty}0.\cr}
$$
Thus \rife{ec1} implies
$$
\eqalign{ &
\la\into w_n\,[(w_n-l)^+]^{2\bar q-1}\,dx+(2\bar q-1)\into
 A(x)\nabla w_n\nabla w_n\,[(w_n-l)^+]^{2\bar q-2}
\,dx \cr &\qquad \leq \frac{C_\de}{\vep^{q-1}} \into
\left|\nabla w_n\right|^q\, [(w_n-l)^+]^{2\bar q-1}\,dx+o(1)_n\,,\cr}
$$
where $o(1)_n$ goes to zero as $n$ tends to infinity.
Neglecting the zero order term which is
positive, and using that $A(x)\geq \alpha I$, we have
\be\label{ast}
\alpha(2\bar q-1)\into
 |\nabla w_n|^2[(w_n-l)^+]^{2\bar q-2}
\,dx \leq \frac{C_\de}{\vep^{q-1}} \into
\left|\nabla w_n\right|^q\,[(w_n-l)^+]^{2\bar
q-1}\,dx+o(1)_n\,.
\ee 
Young's inequality implies
$$
\eqalign{ &
\frac{C_\de}{\vep^{q-1}} \into
\left|\nabla w_n\right|^q\,[(w_n-l)^+]^{2\bar
q-1}\,dx\leq \frac\alpha2(2\bar q-1)\into
 |\nabla w_n|^2[(w_n-l)^+]^{2\bar q-2}
\,dx\cr &\qquad \qquad+ C_{\vep,\de}\into
[(w_n-l)^+]^{2\bar q+\frac{2(q-1)}{2-q}}\,dx\,,\cr}
$$
hence, using that $2\bar q+\frac{2(q-1)}{2-q}=\bar
q2^*$ we get
$$
\frac\alpha2(2\bar q-1)\into
 |\nabla w_n|^2[(w_n-l)^+]^{2\bar q-2}
\,dx \leq  C_{\vep,\de} \into
 \,[(w_n-l)^+]^{ \bar
q 2^*}\,dx+o(1)_n\,.
$$
Using Sobolev inequality in the left hand side we
obtain
$$
C\left(\into
 \,[(w_n-l)^+]^{ \bar
q 2^*}\,dx\right)^{1-\frac2N}\leq  C_{\vep,\de} \into
 \,[(w_n-l)^+]^{ \bar
q 2^*}\,dx+o(1)_n
$$
We let now $n$ tend to infinity; since $u$, $v$ and 
$\vfi$ all belong to $\elle{\bar q2^*}$, we have
that 
$$
w_n\to w\,:=\, u-(1-\vep)v-\vep \vfi
\quad\hbox{strongly in
$\elle{\bar q2^*}$ as $n$ tends to infinity,}
$$
hence we get
$$
C\left(\into
 \,[(w -l)^+]^{ \bar
q 2^*}\,dx\right)^{1-\frac2N}\leq  C_{\vep,\de} \into
 \,[(w-l)^+]^{ \bar
q 2^*}\,dx\,.
$$
 Since
$1-\frac2N<1$, last inequality implies  that $w\leq
0$; indeed, if $\sup w>0$ (even possibly
infinite), one gets a  contradiction by letting   
$l$ converge to
$\sup w$   and using that
$[(w-l)^+]^{
\bar q 2^*}$ would tend to zero in $\elle1$.

The conclusion is then that 
$w\leq 0$, i.e.
$$
u\leq (1-\vep)v+\vep \vfi\,,
$$
and, letting $\vep\to 0$, $u\leq v$ in $\Omega$. 
\qed
\vskip1em

Let us now deal with the case $\la=0$. Indeed, the same
  proof can be applied, provided there
exists a solution of 
\be\label{auxla}
\cases{-\dive(A(x)\nabla
\vfi)=(\gamma+\de)|\nabla \vfi|^q+f(x)& in
$\Omega$,\cr
\vfi=0&on $\partial \Omega$,\cr}
\ee
for some $\de>0$. This requires a further assumption, which is   a
sort of size condition on the data. 

Indeed, it is  known from \cite{GMP} that there exists a  constant $C_*$, only depending on
$q$ and $N$, such that if 
\be\label{picc}
b^{\frac1{q-1}}\,\|f\|_{\elle{\frac N{q'}}}< \alpha^{q'} \,C_*
\ee
then the problem
\be\label{auxb}
\cases{-\dive(A(x)\nabla
z)=b|\nabla z|^q+f(x)& in
$\Omega$,\cr
z=0& on $\partial \Omega$\cr}
\ee
admits a  solution $z$ such that $(1+|z|)^{\bar q-1}z\in \acca$. In particular, if we
fix $\alpha$ (the coercivity constant of $A(x)$) and $f$, the set 
$$
{\cal B}_f\,:\,=\{b>0\,:\,\hbox{problem \rife{auxb} has a solution $z$\ : $(1+|z|)^{\bar
q-1}z\in \acca$}\}
$$ 
is non empty. Indeed, it is not difficult to see that $  {\cal B}_f$ is even an interval.
In order to assure that \rife{auxla} has a  solution for a certain $\de>0$, we are then
led to  assume that \rife{h2} holds with  $\gamma< \sup \,{\cal B}_f$.

\begin{theo}\label{comp2}
 Let $\la=0$. Assume \rife{matr}, \rife{h1} and  \rife{h2} with
$q\geq 1+\frac 2N$ and $\gamma<\sup \,{\cal B}_f$, which is  defined above.  If
$u$ and
$v$ are respectively a  subsolution and a 
super-solution of \rife{pb} such that $(1+|u|)^{\bar
q-1}u\in \acca$ and $(1+|v|)^{\bar q-1}v\in \acca$, then we
have $u\leq v$  in $\Omega$. In particular  problem
\rife{pb} has a unique solution $u$ such that
$(1+|u|)^{\bar q-1}u \in
\acca$.
\end{theo}

The result of Theorem \ref{comp2} may be rephrased more explicitly in terms of a size
condition on the norm of $f$. Indeed, let $C_*$ be the maximal possible choice in
\rife{picc}, i.e. 
$$
\begin{array}{rl}
C_*=&\sup\{ C>0\,:\, \hbox{if  
$\alpha^{-q'} b^{\frac1{q-1}}\,\|f\|_{\elle{\frac N{q'}}}<C$ then problem \rife{auxb}}\\ &
\qquad\quad \hbox{  has a  solution $z$ : $(1+|z|)^{\bar q-1}z\in \acca$}\}
\end{array}
$$
Then we  have

\begin{corol}\label{coc}
 Let $\la=0$. Assume \rife{matr}, \rife{h1} and  \rife{h2}
with  $q\geq 1+\frac 2N$ and 
\be\label{piccola}
\alpha^{-q'}\gamma^{\frac1{q-1}}\,\|f\|_{\elle{\frac N{q'}}}< C_*
\ee
Then problem \rife{pb} has a  unique solution $u$ such that 
$(1+|u|)^{\bar q-1}u \in
\acca$.
\end{corol}

\vskip1em
\begin{rem}\rm Applying Theorem \ref{comp} and
Corollary  \ref{coc} to the model problem
\begin{equation}\label{pbmod2}
\left\{\begin{array}{l}
\lambda
u-\mbox{div}(A(x)Du)=\gamma\,|Du|^q+f(x)\mbox{ in
}\Omega ,\\ u=0\mbox{ on }\partial\Omega\end{array}
 \right.\end{equation} 
we obtain the results stated in the Introduction.
Observe that if
$\gamma>0$ and
$f\in
\elle{\frac N{q'}}$, one can easily prove that {\it
any} weak solution satisfies $(u^{-})^{\bar q}\in
\acca$; in particular, in that case uniqueness
holds  in the class of solutions $u\in\acca$ such
that $(u^+)^{\bar q}\in \acca$. Clearly, when $\gamma $ is negative we should apply the result to the equation satisfied by $-u$.
\end{rem}

\vskip0,5em
\begin{rem}\rm When considering the case $\la=0$,  a more careful look at the
proof of Theorem \ref{comp} shows that in the inequality  \rife{de} one could  replace
$f$ with
$\tilde f=\sup\limits_{\xi} (H(x,\xi)-\gamma|\xi|^q)^+$. The size condition of Theorem
\ref{comp2} and Corollary \ref{coc} would then concern $\tilde f$ instead of $f$.
If one looks at the model problem \rife{pbmod2}, this simply means that  if $\la=0$
and $\gamma >0$,  the required size condition   only concerns  $f^+$, as it is expected.  
\end{rem}

\vskip0,5em

\begin{rem}\rm When $q\to 2$, the exponent $\bar q\to +\infty$. In fact, if $q=2$
uniqueness for problems like  \rife{pbmod2} holds in the class of solutions $u\in \acca$
such that $e^{\mu u}-1 \in \acca$ for some suitable $\mu>0$. This result is proved (in a 
more general framework) in 
\cite{Bar2}: the idea is to use the change of unknown function $v=e^{\gamma u}-1$, so
that the  standard choice is to take $\mu=\gamma$ and to prove uniqueness when
$v\in \acca$. Otherwise one should
take
$\mu=n\gamma$ for some
$n>1$; in that case one proves uniqueness for solutions such that $|v|^{n-1}v\in \acca$.
However, we point out that
this requires to apply to the equation of $v$  a similar truncation argument as in
Lemma \ref{ren}.
\end{rem}
\vskip1em

\subsection{Comments and  extensions}

\begin{itemize}

\item[\bf 1.] {\bf Data in $W^{-1,r}$.}

The results of this section still hold if the right hand side in \rife{pb} is replaced by
$H(x,\nabla u)+ \dive (g(x))$ with $g(x)\in \elle{N(q-1)}$. 

\item[\bf 2.] {\bf Neumann boundary conditions.}

Our method easily extends to prove  a  comparison principle 
for the homogeneous Neumann problem which can be written in a strong form as
\be\label{neu}
\cases{\la u -\dive(A(x)\nabla u)= H(x,\nabla u)& in $\Omega$, \cr A(x)\nabla u\cdot  \nu(x)=0 &on $\partial \Omega$,\cr}
\ee
where $\nu (x)$ is the outward, unit normal vector to $\partial \Omega$ at $x$.
Of course, we use the classical weak formulation which says that $u\in H^1(\Omega)$ is  a weak solution of \rife{neu} if 
$$
\la\into  u \vfi\,dx+\into A(x)\nabla u\nabla \vfi\,dx=\into  H(x,\nabla u)\vfi\,dx\quad
\forall \vfi\in H^1(\Omega)\cap\elle\infty.
$$
Then one has 
\begin{theo} \label{neuu}
Assume \rife{matr}, \rife{h1}, \rife{h2} and
that $\la>0$. Let $q\geq 1+\frac 2N$. If $u$ and $v$
are respectively a  subsolution and a 
super-solution of \rife{neu} such that $(1+|u|)^{\bar
q-1}u\in H^1(\Omega)$ and $(1+|v|)^{\bar q-1}v\in H^1(\Omega)$, then we
have $u\leq v$  in $\Omega$. 

 In particular,  problem
\rife{neu} has a unique weak solution $u$ such that
$(1+|u|)^{\bar q-1}u \in
H^1(\Omega)$.
\end{theo}

\proof The proof follows the same steps as for Theorem \ref{comp}. Note that Lemma
\ref{ren}   is still true without any modification. Then one defines $\vfi$ as a 
solution of 
$$
\cases{\la\,\vfi-\dive(A(x)\nabla
\vfi)=(\gamma+\de)|\nabla \vfi|^q+f(x)& in
$\Omega$,\cr
A(x)\nabla \vfi \cdot  \nu(x)  =0& on $\partial \Omega$,\cr}
$$
and, setting $w_n= T_n(u)-(1-\vep )T_n(v)-\vep \vfi$, one obtains
$$
\eqalign{&
\la\into w_n[(w_n-l)^+]^{2\bar q-1}dx+ \frac\alpha2(2\bar q-1)\into
 |\nabla w_n|^2[(w_n-l)^+]^{2\bar q-2}
\,dx \cr &
\qquad \qquad \qquad \leq  C_{\vep,\de} \into
 \,[(w_n-l)^+]^{ \bar
q 2^*}\,dx+o(1)_n\,.\cr}
$$
Since $\la>0$ one deduces
$$
\| [(w_n-l)^+]^{\bar q}\|_{H^1(\Omega)}\leq \tilde C_{\vep,\de} \into
 \,[(w_n-l)^+]^{ \bar
q 2^*}\,dx+o(1)_n\,.
$$
Using now Sobolev inequality one concludes as in the Dirichlet case.

We only need to require here that $\Omega$  has enough regularity so that the Sobolev
inequality holds. 
\qed

\item[\bf 3.] {\bf Unbounded domains.}

A slight refinement of our proof gives a similar result in case of unbounded domains. To be more precise, let $\Omega$ be  a general domain, not necessarily bounded.
Let still $q\geq 1+\frac2N$ and $\bar q=\frac{(N-2)(q-1)}{2(2-q)} $.

By a  solution of \rife{pb} we mean a function $u$ such that
\be\label{regu}
u\,\psi \in H^1_0(\Omega)\quad \forall \psi\in C_c^\infty(\rn)\,, \quad |u|^{\bar
q-1}u\in \acca
\ee
and
\be\label{eqill}
\la\into u\xi\,dx+ \into A(x)\nabla u\nabla \xi\,dx=\into H(x,\nabla u)\xi\,dx
\quad \forall \xi\in C_c^\infty(\Omega)\,.
\ee
Note that condition \rife{regu} gives a meaning to the Dirichlet condition on
$\partial
\Omega$; roughly speaking, one has  (in a weak sense) $u=0$ on $\partial \Omega\cap B_R$
for any ball
$B_R$ and
$u=0$ at infinity as well (since $\bar q\geq 1$).

The existence of a solution of \rife{pb} in the sense of \rife{regu}--\rife{eqill} has
been proved in \cite{GMP}. It was also pointed  out that, due to the regularity of
$|u|^{\bar q-1}u$, one can allow in \rife{eqill} any test function $\xi$ of the form
$S(u)$, where $S(0)=0$ and $|S'(t)|\leq |t|^{2\bar q-2}$. This is achieved by choosing
$\xi=S(u)\,\zeta(\frac{|x|}n)$, where $\zeta\in C_c^\infty(B_2)$, $\zeta\equiv 1$ on
$B_1$, and  letting $n$ go to infinity, which is allowed  thanks to  \rife{regu} and
\rife{h1}--\rife{h2}.

\begin{theo}  Assume \rife{matr}, \rife{h1}, \rife{h2} with $q\geq 1+\frac
2N$, and that

(i) either $\la >0$

(ii) or $\la=0$ and  
\rife{piccola} holds true.

Then there exists a unique solution $u$ of \rife{pb} in the sense of
\rife{regu}--\rife{eqill}.
\end{theo}

\proof Note that Lemma \ref{ren} still holds true, i.e. \rife{req} holds for any $\xi\in
C_c^\infty(\rn)$, and estimate \rife{tail} is still valid, since it only depends on the
fact that
$|u|^{\bar q-1}u\in
\acca$.  We proceed then as in Theorem \ref{comp}: let $\vfi$ be a  solution (whose
existence is proved in \cite{GMP}) of the auxiliary problem
$$
\cases{\la\,\vfi -\dive(A(x)\nabla
\vfi )=(\gamma+\de)|\nabla \vfi |^q+ f(x)& in
$\Omega$,\cr
|\vfi|^{\bar q-1}\vfi\in \acca\,,\quad \hbox{$\vfi\psi\in \acca$ for any $\psi\in
C_c^\infty(\rn)$.}& \cr}
$$
Defining 
$w_n= T_n(u)-(1-\vep )T_n(v)-\vep \vfi$ we obtain
$$
\eqalign{  
\la\into w_n\,\xi \,dx+\into
 A(x)\nabla w_n\nabla \xi 
\,dx \leq &\frac{C_\de}{\vep^{q-1}} \into
\left|\nabla w_n\right|^q\,\xi \,dx\cr &\qquad  +\langle
I_{n}^{u},\xi\rangle-(1-\vep)\langle I_{n}^{v},\xi\rangle\,,\cr}
$$
for any $\xi\in C_c^\infty(\Omega)$. By density, one can allow $\xi=z\,\psi$ for any
$\psi\in C_c^\infty(\rn)$ and for any $z\in \elle\infty$ such that $z\psi\in \acca$. 

Now choose
$\xi=[(w_n-l)^+]^{2\bar q-1}\psi_j^2$, where $\psi_j=\psi\left(\frac{|x|}j\right)$,  
$\psi\in C_c^\infty(B_2)$, $\psi\equiv 1$ on $B_1$. We get
$$
\eqalign{  &
\la\into w_n\,[(w_n-l)^+]^{2\bar
q-1}\psi_j^2 \,dx+\alpha(2\bar q-1)\into
 |\nabla w_n|^2[(w_n-l)^+]^{2\bar
q-2}\psi_j^2
\,dx\cr & \leq \frac{C_\de}{\vep^{q-1}} \into
\left|\nabla w_n\right|^q\,[(w_n-l)^+]^{2\bar
q-1}\psi_j^2\,dx
\cr &\qquad    -
2\into A(x)\nabla w_n\nabla \psi_j\,[(w_n-l)^+]^{2\bar
q-1}\psi_j\,dx \cr  &\qquad \quad+\langle I_{n}^{u},[(w_n-l)^+]^{2\bar
q-1}\psi_j^2\rangle-(1-\vep)\langle
I_{n}^{v},[(w_n-l)^+]^{2\bar
q-1}\psi_j^2\rangle\,,\cr}
$$
Since $\psi_j\leq 1$, and due to estimate \rife{tail} we have
$$
\eqalign{ &
\langle I_{n}^{u},[(w_n-l)^+]^{2\bar
q-1}\psi_j^2\rangle-(1-\vep)\langle
I_{n}^{v},[(w_n-l)^+]^{2\bar
q-1}\psi_j^2\rangle\cr &\qquad \leq 
n^{2\bar q-1}[\|I_n^{u}\|_{\elle1}+\|I_n^{v}\|_{\elle1}]=o(1)_n\cr}
$$
Using Young's inequality we get
\be\label{prij}
\eqalign{  &
\la\into w_n\,[(w_n-l)^+]^{2\bar
q-1}\psi_j^2 \,dx+\frac{\alpha(2\bar q-1)}2\into
 |\nabla w_n|^2[(w_n-l)^+]^{2\bar
q-2}\psi_j^2
\,dx\cr & \leq   C_{\de,\vep} \into
[(w_n-l)^+]^{2^*\bar
q}\psi_j^2\,dx
 + C
\into |\nabla \psi_j|^{2}\,[(w_n-l)^+]^{2\bar
q}\,dx +o(1)_n\,.\cr}
\ee
Observe that $w_n$ belongs to $\elle{2^*\bar q}$, since it is so for $u$, $v$ and
$\vfi$. Moreover $|\nabla
\psi_j|^{2}$ weakly converges to zero in $\elle{\frac N2}$, so that
$$
\lim\limits_{j\to +\infty}\into |\nabla \psi_j|^{2}\,[(w_n-l)^+]^{2\bar
q}\,dx=0\,.
$$
Since \rife{prij} implies
$$
\eqalign{&\qquad 
\| [(w_n-l)^+]^{\bar q}\psi_j\|_{\acca}^2 \leq C_{\de,\vep} \into
[(w_n-l)^+]^{2^*\bar
q}\psi_j^2\,dx
\cr & + C
\into |\nabla \psi_j|^{2}\,[(w_n-l)^+]^{2\bar
q}\,dx +o(1)_n\,,\cr}
$$
passing to the limit   as $j$ goes to infinity we find then that
$[(w_n-l)^+]^{\bar q}\in
\acca$  and
$$
\| [(w_n-l)^+]^{\bar q}\|_{\acca}^2\leq C_{\de,\vep} \into
[(w_n-l)^+]^{2^*\bar
q} \,dx
 + o(1)_n\,. 
$$
Using Sobolev inequality and that
$$
w_n\to   u-(1-\vep)v-\vep \vfi
\quad\hbox{strongly in
$\elle{\bar q2^*}$ as $n$ tends to infinity,}
$$
letting $n$ go to infinity the conclusion follows as
in Theorem \ref{comp}.
\qed

Finally, when $\la>0$ a  similar result can  be given in case of Neumann boundary
conditions proceeding as in Theorem \ref{neuu}.

\end{itemize}

\vskip1em
\section{The case ${\bf q\leq 1+\frac2N}$.}

\vskip1em

We start by extending  Theorem \ref{comp} to the  case $q<1+\frac 2N$. 
However, in view of \rife{h2} and \rife{ac}, in this case solutions do not belong
in general to
$\acca$, so that one needs first to define a 
suitable concept of solution. It seems useful to
adopt the notion  of renormalized solutions; this notion, introduced first in  
\cite{DL} for transport equations, has been adapted   to second order elliptic
equations in  \cite{BDGM}, \cite{LM}, and recently used in several other contexts when
dealing with unbounded solutions having infinite energy.

\noindent   Let us recall that the auxiliary functions $T_n(s)$ are defined in
\rife{deft}.

\begin{defi} A  renormalized solution of problem \rife{pb}
is a function $u\in \elle1$ such that   $T_n(u)\in\acca$ for any $n>0$, $H(x,\nabla u)\in
\elle1$ and which satisfies
\be\label{RE}
\eqalign{ &
\la\into u\,S(u)\xi\,dx+\into
A(x)\nabla u\nabla (S(u)\xi)\,dx  =\into
H(x,\nabla  u)S(u) \xi\,dx \cr}
\ee
for any Lipschitz function $S$ having compact support
and  for any $\xi \in H^1(\Omega)\cap \elle\infty$ such that  $S(u)\xi\in\acca$.

Renormalized sub or super-solutions are defined in the same way by replacing the equality
in \rife{RE} with the suitable  inequality.
\end{defi}

Clearly, if $u\in \acca$ is a  weak solution then it is also a renormalized
solution: indeed, one can choose  $S(u)\xi\in \acca\cap \elle\infty$   as test function in
\rife{WE} and obtain \rife{RE}. Thus, for $\acca$ solutions, the weak and renormalized
formulations are equivalent. 
However, as in the previous section,  we deal with  
solutions
$u$ such that
$(1+|u|)^{\bar q-1}u\in \acca$, where $\bar q=\frac{(N-2)(q-1)}{2(2-q)} $: if $q<
1+\frac2N$ then  $\bar q<1$, so that solutions do not have finite
energy (i.e. they are not in $\acca$). In this case, the renormalized formulation
 is meant to allow 
test functions depending on $u$ itself, which can not be ensured by using the
simpler
 distributional
formulation. Another possible formulation based on a  duality argument is mentioned later
(see \rife{dual}).

The existence of a renormalized solution $u$ such that
$(1+|u|)^{\bar q-1}u\in \acca$ has been proved in
\cite{GMP}. 
The method of proof given in Section 2 can be easily adapted to provide uniqueness of
such solutions.

\begin{theo}\label{compre}
Assume \rife{matr}, \rife{h1}, \rife{h2} with $q<1+\frac 2N$. 
Let $\la>0$. If $u$ and $v$ are respectively   renormalized subsolution
and    super-solution of \rife{pb} such that $(1+|u|)^{\bar
q-1}u\in \acca$ and $(1+|v|)^{\bar q-1}v\in \acca$, then we
have $u\leq v$ in $\Omega$. 

In particular,  problem
\rife{pb} has a unique renormalized solution $u$ such that
$(1+|u|)^{\bar q-1}u \in
\acca$.
\end{theo}

\proof  First we observe that Lemma \ref{ren} still holds for renormalized solutions:
indeed, choosing  in \rife{RE} $S=\theta_n$ (see \rife{deft}) yields
the same as  \rife{renp}, so that we have
\be\label{rn} 
\eqalign{ &
\la\into
T_n(u)\xi\,dx+\into A(x)\nabla
T_n(u)\nabla \xi \,dx\cr  &\leq \into H(x,\nabla 
T_n(u)) \xi\,dx+
\langle I_{n},\xi\rangle\,,\cr}
\ee
where $I_n$ is defined as in \rife{mun}. Moreover, proceeding exactly as in 
 Lemma \ref{ren} we obtain the estimate  
\be\label{tail2}
\eqalign{&
\lim\limits_{n\to +\infty} n^{2\bar q-1}\,\|
I_{n}\|_{L^1(\Omega)}=0\,.
\cr}
\ee
The same can be proved as regards $v$. Then, using 
the convexity of
$H$, we can proceed as in the proof of Theorem \ref{comp}, in order to obtain that
 \be\label{ecre}
\eqalign{ &
\la\into [T_n(u)-(1-\vep)T_n(v)-\vep \vfi]\,\xi \,dx\cr &\quad +\into
 A(x)\nabla [T_n(u)-(1-\vep)T_n(v)-\vep \vfi]\nabla \xi 
\,dx \cr&\qquad \leq \frac{C_\de}{\vep^{q-1}} \into
\left|\nabla [T_n(u)-(1-\vep)T_n(v)-\vep \vfi]\right|^q\,\xi \,dx\cr
 &\quad\quad
-\vep
\left[\into \la\vfi\,\xi +A(x)\nabla
\vfi\nabla
\xi\,dx-(\gamma+\de)|\nabla
\vfi|^q\xi-f(x)\xi\right]\cr 
& \qquad \qquad +\langle I_{n}^{u},\xi\rangle-(1-\vep)\langle I_{n}^{v},\xi\rangle\,,\cr}
\ee
for any $\vfi$,  $\xi\in \acca\cap\elle\infty
$, $\xi\geq 0$. We
define here $\vfi_n $ to be  a solution of  
\be\label{auxn}
\cases{\la\,\vfi_n-\dive(A(x)\nabla
\vfi_n)=(\gamma+\de)|\nabla \vfi_n|^q+T_n(f(x))& in
$\Omega$,\cr
\vfi_n=0&on $\partial \Omega$.\cr}
\ee
Note  that $\vfi_n$ is nonnegative and  belongs to $\acca\cap
\elle\infty$.  It is proved in \cite{GMP} that $(1+\vfi_n)^{\bar q-1}\vfi_n$ is
bounded in
$\acca$ and 
\be\label{fin}
\vfi_n\to \vfi\qquad \hbox{strongly in $\elle{\bar q\,2^*}$,}
\ee
where $\vfi$ is a renormalized solution of \rife{auxn} corresponding to $f$, and
satisfying $(1+\vfi)^{\bar q-1}\vfi\in \acca$.

Setting
$$
w_n=T_n(u)-(1-\vep)T_n(v)-\vep \vfi_n
$$
and using the equation satisfied by $\vfi_n$ we obtain from \rife{ecre}
\be\label{prewn}
\eqalign{ &
\la\into w_n\,\xi \,dx+\into
 A(x)\nabla w_n\nabla \xi 
\,dx \leq \frac{C_\de}{\vep^{q-1}} \into
\left|\nabla w_n\right|^q\,\xi \,dx\cr
 &\quad\quad
-\vep
\into (T_n(f)-f)\xi\,dx+\langle I_{n}^{u},\xi\rangle-(1-\vep)\langle
I_{n}^{v},\xi\rangle\,.\cr}
\ee
Note
that since $\frac N{N-1}<q< 1+\frac 2N$ then the
exponent $\bar q\in (\frac12, 1)$, hence $2\bar
q-1\in (0, 1)$; for this reason we choose now $\xi=\xi_{n,\sigma}$ with
$$
\xi_{n,\sigma}= [\sigma +
(w_n-l)^+]^{2\bar q-1}-\sigma^{2\bar q-1},
$$
where $l$, $\sigma>0$.  We have, using \rife{matr},
that
$$
\liminf\limits_{\sigma \to 0}\into
 A(x)\nabla w_n\nabla\xi_{n,\sigma}
\,dx  \geq \frac{(2\bar q-1)\alpha}{\bar q^2}\into |\nabla [(w_n-l)^+]^{\bar q}|^2\,dx\,.
$$
Note that $\xi_{n,\sigma}\leq [(w_n-l)^+]^{2\bar q-1}$, and clearly $\xi_{n,\sigma}\to
[(w_n-l)^+]^{2\bar q-1}$ as $\sigma \to 0$. From \rife{prewn} we obtain, as $\sigma \to0$,
$$
\eqalign{ &
\la\into w_n\,[(w_n-l)^+]^{2\bar q-1}\,dx+\frac{(2\bar q-1)\alpha}{\bar q^2}\into |\nabla
[(w_n-l)^+]^{\bar q}|^2\,dx \cr &\quad   
\leq
\frac{C_\de}{\vep^{q-1}} \into
\left|\nabla w_n\right|^q\,[(w_n-l)^+]^{2\bar q-1}\,dx
\cr &\qquad -\vep
\into (T_n(f)-f)[(w_n-l)^+]^{2\bar q-1}\,dx\cr
 &\qquad\qquad+\langle I_{n}^{u},[(w_n-l)^+]^{2\bar
q-1}\rangle-(1-\vep)\langle I_{n}^{v},[(w_n-l)^+]^{2\bar q-1}\rangle\,.\cr}
$$
Since $(w_n-l)^+\leq 2n$, using \rife{tail2} we obtain 
that last two terms go to zero as $n$ tends to infinity.
Moreover, since $u$ and $v$ belong to $\elle{\bar q\,2^*}$ and using \rife{fin}, we have
that 
$ [(w_n-l)^+]^{2\bar q-1}$ converges strongly in $\elle{\frac {\bar q\,2^*}{2\bar q-1}}$;
but we have $\frac {\bar q\,2^*}{2\bar q-1}= (\frac N{q'})'$, and since $T_n(f)-f$
strongly converges to zero in $\elle{\frac N{q'}}$, we conclude that
$$
\eqalign{ &
\la\into w_n\,[(w_n-l)^+]^{2\bar q-1}\,dx+\frac{(2\bar q-1)\alpha}{\bar q^2}\into |\nabla
[(w_n-l)^+]^{\bar q}|^2\,dx \cr &\qquad   
\leq
\frac{C_\de}{\vep^{q-1}} \into
\left|\nabla w_n\right|^q\,[(w_n-l)^+]^{2\bar q-1}\,dx
+o(1)_n\,,\cr}
$$
where $o(1)_n$ goes to zero as $n$ tends to infinity.  This inequality is  the same  as  \rife{ast}, and the conclusion of the proof  is
exactly as in Theorem \ref{comp}.
\qed

A  similar result holds in case $\la=0$ if the data satisfy a suitable size condition,
following the same principle as in Theorem
\ref{comp2}. We leave the details to the reader.
\vskip1em
We are going now to see a different approach to uniqueness, which is based on a simpler
linearization principle. This  approach, which was not possible in the situation of
Section 2, is allowed if $q\leq 1+\frac2N$ (note that the limiting value $q=1+\frac 2N$ is
included too), and provides uniqueness in   a more general context. Namely, we consider
the problem
\be\label{pbnl}
\cases{\la u-\dive(a(x,\nabla u))= H(x,\nabla u)& in $\Omega$,\cr
u=0&on $\partial
\Omega$\cr}
\ee
where $a(x,\xi)\,:\,\Omega\times\rn\to \rn$ is a Carath\'eodory function such
that 
\be\label{au}
[a(x,\xi)-a(x,\eta)]\cdot (\xi-\eta) \geq \alpha |\xi-\eta|^2
\ee
\be\label{ad}
|a(x,\xi)|\leq \beta (k(x)+|\xi|)\,,\quad \beta>0\,, \quad k(x)\in \elle2\,.
\ee
We assume that  $H(x,\xi)\,:\,\Omega\times \rn\to \r$ is a Carath\'eodory function 
which satisfies 
\be\label{hlip}
\eqalign{ &
| H(x,\xi)-H(x,\eta)|\leq \gamma( b(x)+ |\xi|^{q-1}+|\eta|^{q-1})\,|\xi-\eta|\,,
\cr & \qquad \qquad\qquad \quad
b(x)\in \elle N\,,\,\gamma>0,\cr}
\ee
and
\be\label{h0}
H(x,0)\in \elle {\frac N{q'}}.
\ee
Note that assumptions \rife{hlip} and \rife{h0} imply that $H(x,\xi)$ satisfies the growth
condition
\rife{h2}. On the other hand, no convexity is now assumed on $H(x,\cdot)$.

As in Definition \ref{RE}, we say that  a function $u\in \elle1$ is a  renormalized
subsolution (super-solution) of problem \rife{pbnl} if  
$T_n(u)\in\acca$ for any $n>0$, $H(x,\nabla u)\in \elle1$ and  
\be\label{REnl}
\eqalign{ &
\la\into u\,S(u)\xi\,dx+\into
a(x,\nabla u)\nabla (S(u)\xi)\,dx \leq\, (\geq) \,\,\into
H(x,\nabla  u)S(u) \xi\,dx \cr}
\ee
for any Lipschitz function $S$ having compact support
and for any $\xi \in H^1(\Omega)\cap \elle\infty$ such that   $S(u)\xi\in\acca$.

\vskip1em

We start with two important properties of solutions in the class \rife{clac}. We will
need a  slight modification of the truncation functions $T_n(s)$. Namely, we set
$$
{\cal T}_n(t)=\int_0^t \ct_n'(s)ds\,,\quad \ct_n'(s)=\cases{1&if $|s|<n$\cr
n+1-|s|  &if $n<|s|<n+1$,\cr 0&if
$|s|>n+1$\cr}
$$

\begin{lem}\label{renno}  Assume \rife{au}--\rife{h0} with $\frac N{N-1}<q\leq 1+\frac 2N$.
Let 
$u$ be a renormalized  subsolution  of \rife{pbnl} such that $(1+|u|)^{\bar
q-1}u\in \acca$, with $\bar
q=\frac{(N-2)(q-1)}{2(2-q)}$. Then we have
\begin{itemize}
\item[(i)]
\be\label{bg}
u\in \sob {N(q-1)}.
\ee
\item[(ii)] 
\be\label{old}
\lim\limits_{n\to +\infty}\quad 
n^{2\bar q-1}\!\!\!\!\!\int\limits_{\{n<|u|<n+1\}}a(x,\nabla u)\nabla
u\,dx =0
\ee
\item[(iii)] 
for every $\xi\in \acca\cap
\elle\infty$ and for every $n$
\be\label{renno1}
\eqalign{ &
\la\into
\ct_n(u)\xi\,dx+\into a(x,\nabla
u)\ct_n'(u)\nabla \xi \,dx\cr &\quad \leq \into H(x,\nabla 
\ct_n(u)) \xi\,dx+
\langle I_{n},\xi\rangle\,,\cr}
\ee
with
\be\label{tail3}
\eqalign{&
\lim\limits_{n\to +\infty} n^{2\bar q-1}\,\|
I_{n}\|_{L^1(\Omega)}=0\,.\cr}
\ee
\end{itemize}
\end{lem}

\proof
The regularity \rife{bg} follows directly from the fact that $(1+|u|)^{\bar
q-1}u\in \acca$.  This was first observed, in  a different context, in \cite{BG}; for the
reader's convenience, we  recall the simple argument. Indeed,   due to Sobolev and
H\"older 's inequalities, we have
$$
\eqalign{ &
\left(\into (|u|)^{(N(q-1))^*}\,dx \right)^{2-q}
\leq 
\into |\nabla
u|^{N(q-1)}dx\cr &\quad \leq \left(\into  \frac {|\nabla u|^2}{(1+|u|)^{2-2\bar
q}}\right)^{\frac {N(q-1)}2} 
\left( \into (1+|u|)^{\frac{2N(q-1)(1-\bar q)}{2-N(q-1)}}dx \right)^{1-\frac{N(q-1)}2}\cr}
$$
Since, by definition of $\bar q$, we have  $\frac{2N(q-1)(1-\bar
q)}{2-N(q-1)}=(N(q-1))^*$ and since $2-q>1-\frac{N(q-1)}2$, we conclude that
$$
\| u\|_{\sob {N(q-1)}} \leq  c (1+ \| (1+|u|)^{\bar
q-1}u\|_{\acca}^{\frac1{\bar q}})\,.
$$
To prove (ii), 
take in \rife{REnl} $\xi=1$ and $S(t)=\theta_{n}(t)\int_0^t|s|^{2\bar
q-1}\chi_{\{n-1<|s|<n\}}ds$, where
$\theta_n$ is defined in \rife{deft}. Since
$S(t)\leq (1+|t|)^{2\bar q-1}\chi_{\{n-1<|u|\}}$ we have
\be\label{pold}
\eqalign{ &
\int\limits_{\{n-1<|u|<n\}} a(x,\nabla u)\nabla u |u|^{2\bar q-1}\,dx \leq 
c\int\limits_{\{n-1<|u|\}} |H(x,\nabla u)|\,\,|u|^{2\bar q-1}\,dx\cr &\qquad 
+\frac cn\int\limits_{\{n<|u|<2n\}}a(x,\nabla u)\nabla u\,|u|^{2\bar q-1}\cr}
\ee
Observe that, by \rife{ad}, 
$$
 \frac1n\int\limits_{\{n<|u|<2n\}}a(x,\nabla u)\nabla u \,|u|^{2\bar q-1}dx\leq 
 c\int\limits_{\{n<|u|<2n\}}[k(x)+|\nabla u|]|\nabla u|\,|u|^{2\bar q-2}dx
$$
which yields, since $\bar q\leq 1$,
$$
\frac1n\int\limits_{\{n<|u|<2n\}}a(x,\nabla u)\nabla u\,|u|^{2\bar q-1}dx\leq 
c\int\limits_{\{n<|u|<2n\}}[k(x)^2+|\nabla (|u|^{\bar q-1}u)|^2]dx
$$
Thus
$$
\lim\limits_{n\to +\infty} 
\frac1n\int\limits_{\{n<|u|<2n\}}a(x,\nabla u)\nabla u\,|u|^{2\bar q-1}=0
$$
Moreover, since $H$ still satisfies \rife{h2}, we have, as in the proof of Lemma
\ref{ren},
$$
\int\limits_{\{n-1<|u|\}}\!\!\!\!\!\!\!\!  |H(x,\nabla u)|\,\,|u|^{2\bar q-1}\,dx\leq
c\int\limits_{\{n-1<|u|\}} \!\!\!\!\!\!\!\! [|\nabla (|u|^{\bar q-1}u)|^2+ (|u|^{\bar
q})^{2^*}]\,dx\quad \mathop{\to}^{n\to +\infty} 0
$$
so that we conclude from \rife{pold}
$$
\int\limits_{\{n-1<|u|<n\}} \!\!\!\!\!\!\!\! a(x,\nabla u)\nabla u\,  |u|^{2\bar q-1}
\quad \mathop{\to}^{n\to +\infty} 0\,
$$
hence \rife{old}.

The proof of (iii) follows the outlines of Lemma \ref{ren}.
Choose $S=\ct'_n(t)$ in \rife{REnl}, so that
$$
\eqalign{ &
\la\into
\ct_n(u)\xi\,dx+\into a(x,\nabla u)
\ct_n'(u)\nabla \xi \,dx  \leq \into H(x,\nabla 
\ct_n(u)) \xi\,dx+
\langle I_{n},\xi\rangle\cr}
$$
where 
 $I_{n}$ is  defined as
\be\label{munbi}
\eqalign{ &
\langle
I_{ n},\xi\rangle=
 \int\limits_{\{n<u<n+1\}}\!\!\!\!
a(x,\nabla u)\nabla u\,\xi\,dx -
\int\limits_{\{-n-1<u<-n\}}\!\!\!\!a(x,\nabla
u)\nabla u\,\xi\,dx \cr &
+\la\into (\ct_n(u)-u\ct'_n(u))\xi\,dx+\into
(H(x,\nabla u)\ct'_n(u)-H(x,\nabla
\ct_n(u)))\xi\,dx\,.\cr}
\ee
As in Lemma \ref{ren}, using the growth condition on $H$ we obtain estimates like
\rife{in3} and
\rife{in4}; moreover, for the first two terms of \rife{munbi} we use \rife{old}. Finally, 
we can conclude that 
\rife{tail3} holds.
\qed

Note that the borderline value $q=1+\frac 2N$ is included in the previous lemma as well
as in the following comparison result. However, some statements would read simpler for
this case: in fact, if $q=1+\frac 2N$ then $\bar q=1$, hence $(1+|u|)^{\bar
q-1}u=u$, which belongs to $\acca$ (and \rife{bg} says the same); in particular, in this
case renormalized solutions are also standard $\acca$ weak solutions.

\vskip1em

\begin{theo}\label{nonlin}
 Assume \rife{au}--\rife{h0} with $\frac N{N-1}<q\leq  1+\frac
2N$. Let 
  $\la\geq 0$. If $u$ and $v$
are respectively a renormalized subsolution and   
super-solution of \rife{pbnl} such that $(1+|u|)^{\bar
q-1}u\in \acca$ and $(1+|v|)^{\bar q-1}v\in \acca$, then we
have $u\leq v$ in $\Omega$. 

In particular,  problem
\rife{pbnl} has a unique renormalized solution $u$ such that
$(1+|u|)^{\bar q-1}u \in
\acca$.
\end{theo}

\proof  From Lemma \ref{renno} we have that
\be\label{unu}
\eqalign{ &
\la\into \ct_n(u)\,\xi dx \into a(x,\nabla \ct_n(u))\nabla \xi dx\leq \into H(x,\nabla
\ct_n(u))\xi dx\cr & 
+\into [a(x,\nabla \ct_n(u))-a(x,\nabla u)\ct_n'(u)]\nabla \xi\,dx
+\langle I_n^{u},\xi\rangle\,,\cr}
\ee
for any $\xi\in \acca\cap \elle\infty$, $\xi\geq 0$.

Similarly we deal with the equation satisfied by $v$, so that
\be\label{unv}
\eqalign{ &
\la\into \ct_n(v)\,\xi dx \into a(x,\nabla \ct_n(v))\nabla \xi dx\geq \into H(x,\nabla
\ct_n(v))\xi dx\cr &
+ \into [a(x,\nabla \ct_n(v))-a(x,\nabla v)\ct_n'(v)]\nabla \xi\,dx +\langle
I_n^{v},\xi\rangle\,\cr}
\ee
where 
\be\label{tail4}
n^{2\bar q-1}\|I_n^{v}\|_{\elle1} \mathop{\to}^{n\to +\infty} 0.
\ee
For $k>0$, let us set $G_k(s)=(s-k)^+$: subtracting \rife{unv} from \rife{unu} and choosing
$$
\xi=[G_k(\ct_n(u)-\ct_n(v)) +\sigma]^{2\bar q-1}-\sigma^{2\bar q-1}\,,\quad\sigma>0
$$
 we get
\be\label{uu}
\eqalign{ &
\la\into (\ct_n(u)-\ct_n(v))\,\left([G_k(\ct_n(u)-\ct_n(v))+\sigma]^{2\bar
q-1}-\sigma^{2\bar q-1}\right)dx\cr & 
+ \into [a(x,\nabla \ct_n(u))-a(x,\nabla \ct_n(v))]\nabla
  [G_k(\ct_n(u)-\ct_n(v))+\sigma]^{2\bar q-1} dx\cr &
\leq
\into \![H(x,\nabla \ct_n(u))-H(x,\nabla
\ct_n(v))]\left([G_k(\ct_n(u)-\ct_n(v))+\sigma]^{2\bar q-1}-\sigma^{2\bar
q-1}\right)dx\cr & 
+ \into [a(x,\nabla \ct_n(u))-a(x,\nabla u)\ct_n'(u)]\nabla
  [G_k(\ct_n(u)-\ct_n(v))+\sigma]^{2\bar q-1}\cr & + \into
[a(x,\nabla \ct_n(v))-a(x,\nabla v)\ct_n'(v)]\nabla  
[G_k(\ct_n(u)-\ct_n(v))+\sigma]^{2\bar q-1}
\cr &+\langle  | I_n^{u}|,[G_k(\ct_n(u)-\ct_n(v))+\sigma]^{2\bar q-1}\rangle +\langle
|I_n^{v}|,[G_k(\ct_n(u)-\ct_n(v))+\sigma]^{2\bar q-1}\rangle\cr}
\ee
Since $[G_k(\ct_n(u)-\ct_n(v))]^{2\bar
q-1}\leq c\,n^{2\bar q-1}$ last two terms go to zero as $n$ tends to infinity
thanks to \rife{tail3} and \rife{tail4}. 
Moreover we have from \rife{ad}
$$
\eqalign{ &
|a(x,\nabla \ct_n(u))-a(x,\nabla u)\ct'_n(u)| \cr &\leq [|a(x,\nabla \ct_n(u))|+|a(x,\nabla
u)|]\chi_{\{n<|u|<n+1\}}+ a(x,0)\chi_{\{n+1<|u|\}}\cr &\leq c[|\nabla
u|\chi_{\{n<|u|<n+1\}} + k(x)\chi_{\{n<|u|\}}]\,.
\cr}
$$
Using \rife{au} and Young's inequality, we have
\be\label{uuu}
\eqalign{ &
\la\into (\ct_n(u)-\ct_n(v))\,\left([G_k(\ct_n(u)-\ct_n(v))+\sigma]^{2\bar
q-1}-\sigma^{2\bar q-1}\right)dx\cr & 
+ \into |\nabla G_k(\ct_n(u)-\ct_n(v))|^2
  [G_k(\ct_n(u)-\ct_n(v))+\sigma]^{2\bar q-2} dx\cr &
\leq c
\!\into \![H(x,\nabla \ct_n(u))-H(x,\nabla
\ct_n(v))]\left([G_k(\ct_n(u)-\ct_n(v))+\sigma]^{2\bar
q-1}-\sigma^{2\bar q-1}\right) dx\cr & 
+c \into [|\nabla u|\chi_{\{n<|u|<n+1\}}
+ k(x)\chi_{\{n<|u|\}}]^2
  [G_k(\ct_n(u)-\ct_n(v))+\sigma]^{2\bar q-2}\cr &  + c \into [|\nabla
v|\chi_{\{n<|v|<n+1\}} + k(x)\chi_{\{n<|v|\}}]^2
[G_k(\ct_n(u)-\ct_n(v))+\sigma]^{2\bar q-2}
+o(1)_n\cr }
\ee
Thanks to  \rife{old} in Lemma \ref{renno}, and since $0<2\bar q-1\leq 1$, last two terms go to zero as $n$ tends to infinity (for fixed $\sigma>0$).
Thus, using
also that $\la\geq 0$  and \rife{hlip} we have
$$
\eqalign{ &
\into |\nabla G_k( \ct_n(u)-\ct_n(v))|^2
[G_k(\ct_n(u)-\ct_n(v))+\sigma]^{2\bar q-2}
dx\cr &
\leq
c\int\limits_{E_n}\! [b(x)+|\nabla
\ct_n(u)|^{q-1}+|\nabla
\ct_n(v)|^{q-1}]|\nabla (\ct_n(u)-\ct_n(v))|\times\cr
&\qquad\qquad\qquad \times\,[G_k(\ct_n(u)-\ct_n(v))+\sigma]^{2\bar q-1}dx 
+o(1)_n\,,\cr}
$$
where
$$
E_n=\{x\,:\, \ct_n(u)-\ct_n(v)>k\,,\,|\nabla
(\ct_n(u)-\ct_n(v))|>0\}\,.
$$
Using Young's inequality we get
$$
\eqalign{ &
\into |\nabla G_k( \ct_n(u)-\ct_n(v))|^2
[G_k(\ct_n(u)-\ct_n(v))+\sigma]^{2\bar q-2}
dx\cr &
\leq
c\int\limits_{E_n}\!\!
[b(x)+|\nabla
\ct_n(u)|^{q-1}+|\nabla
\ct_n(v)|^{q-1}]^2\times\cr
&\qquad\qquad\qquad\qquad\times\,[G_k(\ct_n(u)-\ct_n(v))+\sigma]^{2\bar
q}dx+o(1)_n\,.\cr}
$$
Using Sobolev inequality and that $u$, $v\in \sob{N(q-1)}$, we deduce
$$
\eqalign{ &
\left(\into  
\left([G_k(\ct_n(u)-\ct_n(v))+\sigma]^{\bar q}-\sigma^{\bar q}\right)^{2^*}
dx\right)^{\frac2{2^*}}
\leq \cr &\qquad \quad \leq c
\left(\int\limits_{E_n}\!\! [b(x)+|\nabla
u|^{q-1}+|\nabla v|^{q-1}]^Ndx\right)^{\frac2N}\times\cr &\qquad \qquad\qquad
\times\left(\into [G_k(\ct_n(u)-\ct_n(v))+\sigma]^{\bar
q\,2^*}dx\right)^{\frac2{2^*}}+o(1)_n\,.\cr}
$$
Letting $n$ tend to infinity  we obtain
$$
\eqalign{ &
  \left(\into  
\left([G_k(u-v)+\sigma]^{\bar q}-\sigma^{\bar q}\right)^{2^*}
dx\right)^{\frac2{2^*}} 
\leq\cr &\qquad\leq c
\left(\int\limits_{\{u-v>k\,,\,|\nabla (u-v)|>0\}}\!\!\!\!\!\!\!\!\!\!\!\!\!\!\!\!\!\!
[b(x)+|\nabla u|^{q-1}+|\nabla v|^{q-1}]^Ndx\right)^{\frac2N}\times\cr &\qquad\qquad\qquad
\qquad
\times\left(\into [G_k(u-v)+\sigma]^{\bar q\,2^*}dx\right)^{\frac2{2^*}}\,,\cr}
$$
and then, as $\sigma \to 0$,
\be\label{Tru}
\eqalign{ &
\qquad\qquad \left(\into  
\left([G_k(u-v)^+]^{\bar q}\right)^{2^*}
dx\right)^{\frac2{2^*}}\cr &
\leq C
\left(\int\limits_{\{u-v>k\,,\,|\nabla (u-v)|>0\}}\!\!\!\!\!\!\!\!\!\!\!\! [b(x)+|\nabla
u|^{q-1}+|\nabla v|^{q-1}]^Ndx\right)^{\frac2N}\cr &\qquad \qquad \times\left(\into
[G_k(u-v)^+]^{\bar q\,2^*}dx\right)^{\frac2{2^*}}\,.\cr}
\ee
From this inequality one can deduce that $u\leq v$ in $\Omega$. Indeed, we argue by
contradiction. Set  $M=\sup(u-v)$; then, should $M$ be positive, even possibly infinite, 
we have
$$
\lim\limits_{k\to M} \,\,{\rm meas }\{x\,:\,u-v>k\,,\,|\nabla (u-v)|>0\}=0\,,
$$
since either $M=+\infty$ or $|\nabla (u-v)=0|$ a.e. on $(u-v)=M$.
Therefore, using that $u$, $v\in \sob{N(q-1)}$, there exists    $k_0<M$ such that
$$
\left(\int\limits_{\{u-v>k_0\,,\,|\nabla (u-v)|>0\}} \!\!\!\!\!\!\!\!\!\! [b(x)+|\nabla
u|^{q-1}+|\nabla v|^{q-1}]^Ndx\right)^{\frac2N}
<\frac1C
$$
 and then \rife{Tru} implies that  $(u-v)\leq k_0$ almost everywhere, getting a 
contradiction with the fact that $k_0<\sup(u-v)$. We conclude that $u\leq
v$.
\qed

We point out that the previous  theorem extends the uniqueness result which is proved in
\cite{BMMP2} assuming  $H(x,0)\in \dacca$ and for solutions   in $ \acca$. However, the
existence of $\acca$ solutions  can not be proved, nor  it
is expected to hold,  under assumption \rife{h0} with $q <1+\frac2N$, so that,
to be consistent with the existence results (see \cite{GMP}) one actually needs to work
with solutions in the class \rife{clac}.

\subsection{Comments and remarks}

\begin{itemize}
\item[\bf 1.] {\bf   The formulation by duality in the linear case}

Consider problem \rife{pb}, where the second-order operator is linear. Instead of using the notion of renormalized solution,
 a  different formulation  can be given by using
 the linear character of the operator.

\begin{defi}\label{sta} (see \cite{St})
A  function $u\in \elle1$ is a solution of \rife{pb} if $H(x,\nabla u)\in \elle1$ and
\be\label{dual}
\begin{array}{rl}&
\la\into u\,\vfi\,dx- \into u\,\dive(A^*(x)\nabla \vfi)\,dx = \into H(x,\nabla
u)\vfi\,dx\,,\\
\noalign{\medskip}
& \hbox{for every $\vfi\in \acca$: $\dive(A^*(x)\nabla \vfi)\in \elle\infty$,}
\end{array}
\ee
where $A^*(x)$ denotes the adjoint matrix of $A(x)$.
\end{defi}

Note that in Definition \ref{sta} only a  minimal regularity is  asked on $u$, by using
the advantage of linearity to integrate  twice by parts. It is well known (see e.g.
\cite{DMOP}) that, since $H(x,\nabla u)\in \elle1$,  any solution in the sense of
Definition \ref{sta}  also satisfies the renormalized formulation 
\rife{REnl}.   We deduce  then   the following

\begin{theo}
Assume \rife{matr}, \rife{hlip}, \rife{h0} with $\frac N{N-1}<q\leq 1+\frac 2N$. Let $\la\geq 0$. Then
there exists a unique function $u$ which is solution of \rife{pb} in the
sense of Definition
\ref{sta} and such that $(1+|u|)^{\bar q-1}u\in \acca$, with $\bar q=
\frac{(N-2)(q-1)}{2(2-q)}$.
\end{theo}

A  similar result can be given in the convex case (i.e. assuming \rife{h1} and \rife{h2}) for any $q$: $\frac N{N-1}<q<2$, since the results of Theorem \ref{comp}
and Theorem \ref{compre} apply to solutions in the sense \rife{dual} which belong to the
class \rife{clac}.

\item[\bf 2.]{\bf  The case $\bf{q \leq \frac N{N-1}}$ and measure data.}

The question of finding a proper class of solutions
where  uniqueness holds is not relevant if   
$q< \frac N{N-1}$ (note that the counterexample given in  \rife{count} holds only for $q>
\frac N{N-1}$). Indeed, asking only $H(x,\nabla u)\in \elle1$,  the solutions of
\rife{oo} are expected to belong to $\sob r$ for any $r<\frac  N{N-1}$, in particular they
already satisfy \rife{cla}. In fact,   uniqueness results when $q<\frac N{N-1} $ have
already been proved, see e.g.  \cite{BMMP1} for
a  result in a general context including nonlinear operators. 

When $q<\frac N{N-1}$ and in case of linear operators, one can even prove uniqueness if
data are bounded measures, using the formulation \rife{dual} and  a simple duality
argument. This was done in \cite{AP} for the Laplace operator, for completeness we 
sketch the result for the general case.

Let
$H(x,\xi)$ satisfy
\be\label{hun}
\eqalign{ &
| H(x,\xi)-H(x,\eta)|\leq \gamma( b(x)+ |\xi|^{q-1}+|\eta|^{q-1})\,|\xi-\eta|\,,
\cr & \qquad \qquad 
\hbox{with $q< \frac N{N-1}$, $b(x)\in \elle r$ for some $r>N$, $\gamma>0$}\cr}
\ee
and
\be\label{hdu}
H(x,0)\in \elle {1}.
\ee
Let $\mu$ be  a bounded Radon measure in $\Omega$. We say that $u$ is a solution
of 
\be\label{mis}
\cases{\la u-\dive(A(x)\nabla u)= H(x,\nabla u)+ \mu& in $\Omega$, \cr u=0&on $\partial
\Omega$,\cr}
\ee
 if $u\in \elle1$, $H(x,\nabla u)\in \elle1$ and 
\be\label{stam}
\begin{array}{rl}&
\la\into u\,\vfi\,dx- \into u\,\dive(A^*(x)\nabla \vfi)\,dx = \into H(x,\nabla
u)\vfi\,dx+ \into \vfi\,d\mu\,,\\
\noalign{\medskip}
& \hbox{for every $\vfi\in \acca$: $\dive(A^*(x)\nabla \vfi)\in \elle\infty$.}
\end{array}
\ee
Note that such  test functions  $\vfi$ are H\"older continuous by means of De
Giorgi--Nash 's results, hence they can be tested against measures. Then we have

\begin{theo}\label{misu}
Assume \rife{matr}, \rife{hun} and \rife{hdu}, and let $\la\geq 0$. Let $\mu$ be a
bounded Radon measure in $\Omega$. Then there exists  a unique solution $u$ of \rife{mis}.
\end{theo}

\proof Let $u_i$, $i=1,2$,  be two solutions of \rife{mis} in the sense of \rife{stam}.
It is known that $u_i\in \sob r$ for any $r<\frac N{N-1}$. Moreover, if $H_n(x,\xi)$ is a
sequence of bounded  functions such that
$H_n(x,\nabla u_i)$ converges to $H(x,\nabla u_i)$ in $\elle1$,  and if $\mu_n$ is a 
sequence of smooth functions converging to $\mu$ in the weak--$*$ topology of measures,
then the solutions of 
\be\label{misn}
\cases{\la u_{i,n}-\dive(A(x)\nabla u_{i,n})= H_n(x,\nabla u_i)+
\mu_n& in
$\Omega$,\cr
 u_{i,n} =0&on $\partial \Omega$\cr}
\ee
converge to $u_i$ in $\sob r$ for any $r<\frac N{N-1}$.
Moreover, one can choose $H_n$ to be $C^1$ and  still satisfying \rife{hun}, and
converging to $H(x,\xi)$ locally uniformly.

Now, since   $u_{i,n}$ belong to $\acca\cap \elle\infty$, we have
\be\label{tun}
\begin{array}{rl}
&
\into (u_{1,n}-u_{2,n})[\la \vfi-\dive(A^*\nabla
\vfi)]\,dx  
\\
\noalign{\medskip}
& = \into \vfi [\int_0^1 \frac{\partial H_n }{\partial \xi}(x,t\nabla
u_{1,n}+(1-t)\nabla u_{2,n})dt]\nabla (u_{1,n}-u_{2,n})\,dx
\\
\noalign{\medskip}
&+\into [H_n(x, \nabla u_1)-H_n(x,\nabla u_{1,n})]\vfi\,dx
\\
\noalign{\medskip}
&\quad -
 \into [H_n(x, \nabla u_2)-H_n(x,\nabla u_{2,n})]\vfi\,dx
\,.
\end{array}
\ee
Set  
$$
p_n(x)= \int_0^1 \frac{\partial H_n }{\partial \xi}(x,t\nabla
u_{1,n}+(1-t)\nabla u_{2,n})dt
$$
and take $\vfi=\vfi_n$ the solution of
$$
\cases{\la \vfi_n-\dive(A^*(x)\nabla \vfi_n)= -\dive(p_n(x)\,\vfi_n)+T_1(u_1-u_2)& in
$\Omega$,
\cr \vfi_n=0&on
$\partial
\Omega$.\cr}
$$
Since  $ H_n$
satisfies \rife{hun}, and using that $\nabla u_{i,n}$ strongly converge in $\sob r$ for
every $r<\frac N{N-1}$,  we have that $H_n(x,\nabla u_{i,n})$ strongly converges to $H(x,\nabla u_i)$ in
$\elle1$, and there exists $\de>0$ such that $p_n$ strongly
converges in
$\elle{N+\de}^N$. By standard regularity results this implies that $\vfi_n$ is uniformly
bounded (even relatively compact) in
$\elle\infty$, hence last two terms in
\rife{tun} converge to zero. Passing to the limit we get
$$
\into (u_1-u_2)T_1(u_1-u_2)\,dx=0
$$
so that $u_1=u_2$.
\qed
\vskip0,5em

\begin{rem}\rm Note that the case $N=2$  also enters in the previous situation; indeed,
when
$N=2$ the values $\frac N{N-1}$ and $1+\frac 2N$ coincide and $\frac N{N-1}=1+\frac
2N=2$. Thus, in the subcritical case  $q<2$   the main uniqueness result reads as in 
Theorem \ref{misu}, at least for linear operators. For nonlinear operators and with data in
$\elle1$, this case is treated in the results in
\cite{BMMP1}.
\end{rem}

\vskip0,5em
\begin{rem}\rm Finally,  the case  $q=\frac N{N-1}$ is  a critical one;
adapting the counterexample
\rife{count} it is still possible to construct a non trivial solution $u$ of the
homogeneous equation
$$
-\Delta u=|\nabla u|^{\frac N{N-1}}\,,\qquad u\in \sob{\frac N{N-1}}
$$ 
so that looking for a smaller class where uniqueness holds is still necessary.
The radial case suggests   that  uniqueness holds here for solutions  $u$ such that
$|\nabla u|\in L^{\frac N{N-1}}(\log L)^{N-1}$. Indeed, in order to extend Theorem
\ref{nonlin} to $q=\frac N{N-1}$  one should work in the context of Orlicz
spaces, and assumption \rife{h0} should be suitably modified as well, e.g. by asking
$H(x,0)\in L^1(\log L)^{N-1}$.
\end{rem}

\end{itemize}

\vskip2em
\noindent{\bf Guy Barles}
\vspace{3mm}\\
{\it Laboratoire de Math\'ematiques et Physique Th\'eorique}\\
{\it UMR CNRS  6083}
\\
{\it Facult\'e des Sciences et Techniques}
\\{\it Universit\'e de Tours} \\
{\it Parc de Grandmont}\\
{\it  Tours 37200, France}\\
\vspace{3mm}\\
{\bf  Alessio  Porretta}\vspace{3mm}\\
{\it Dipartimento di Matematica}
\\
{\it Universit\`a di Roma Tor 
Vergata}\\
{\it Via della Ricerca Scientifica 1}
\\{\it  00133 Roma, Italia}
\end{document}